\documentclass[12pt]{article}
\usepackage{amsmath,amsfonts,amsthm,amssymb,eucal,amsbsy}
\DeclareMathOperator{\R}{Re}
\DeclareMathOperator{\Dom}{Dom}

\newcommand{\Qp}{\mathbb Q_p}
\newcommand{\DA}{D^\alpha}
\newcommand{\DAN}{D^\alpha_N}
\newcommand{\DD}{\mathcal D(\Qp)}
\newcommand{\DDD}{\mathcal D'(\Qp)}

\newcommand{\FF}{\mathcal F}
\numberwithin{equation}{section}
\begin{document}
\newtheorem{lem}{Lemma}
\newtheorem{teo}{Theorem}
\newtheorem{prop}{Proposition}
\newtheorem*{defin}{Definition}
\pagestyle{plain}
\title{Linear and Nonlinear Heat Equations on a $p$-Adic Ball}
\author{ \textbf{Anatoly N. Kochubei}\\ \footnotesize Institute of Mathematics,\\
\footnotesize National Academy of Sciences of Ukraine,\\
\footnotesize Tereshchenkivska 3, Kiev, 01004 Ukraine,\\
\footnotesize E-mail: kochubei@imath.kiev.ua }

\date{}
\maketitle

\bigskip
\begin{abstract}
We study the Vladimirov fractional differentiation operator $D^\alpha_N$, $\alpha >0, N\in \mathbb Z$, on a $p$-adic ball $B_N=\{ x\in \mathbb Q_p:\ |x|_p\le p^N\}$. To its known interpretations via restriction from a similar operator on $\mathbb Q_p$ and via a certain stochastic process on $B_N$, we add an interpretation as a pseudo-differential operator in terms of the Pontryagin duality on the additive group of $B_N$. We investigate the Green function of $D^\alpha_N$ and a nonlinear equation on $B_N$, an analog the classical porous medium equation.
\end{abstract}
\vspace{2cm}


\newpage

\section{Introduction}

The theory of linear parabolic equations for real- or complex-valued functions on the field $\Qp$ of $p$-adic numbers including the construction of a fundamental solution, investigation of the Cauchy problem, the parametrix method, is well-developed; see, for example, the monographs \cite{K2001,Z2016}. In such equations, the time variable is real and nonnegative while the spatial variables are $p$-adic. There are no differential operators acting on complex-valued functions on $\Qp$, but there is a lot of pseudo-differential operators. A typical example is Vladimirov's fractional differentiation operator $\DA$, $\alpha >0$; see the details below. This operator (as well as its multi-dimensional generalization, the so-called Taibleson operator) is a $p$-adic counterpart of the fractional Laplacian $(-\Delta )^{\alpha/2}$ of real analysis.

Already in real analysis, an interpretation of nonlocal operators on bounded domains is not straightforward; see \cite{BV} for a survey of various possibilities. In the $p$-adic case, Vladimirov (see \cite{VVZ}) defined a version $\DAN$ of the fractional differentiation on a ball $B_N=\{ x\in \Qp:\ |x|_p\le p^N\}$ as follows. One takes a test function on $B_N$, extends it onto $\Qp$ by zero, applies $\DA$, and restricts the resulting function to $B_N$. Then it is possible to consider a closure of the obtained operator, for example, on $L^2(B_N)$.

In \cite{K2001} (Section 4.6), a probabilistic interpretation of this operator was given. Let $\xi_\alpha (t)$ be the Markov process with the generator $\DA$ on $\Qp$. Suppose that $\xi_\alpha (0)\in B_N$ and denote by $\xi_\alpha^{(N)}(t)$ the sum of all jumps of the process $\xi_\alpha (\tau )$, $\tau \in [0,t]$ whose $p$-adic absolute values exceed $p^N$. Consider the process $\eta_\alpha (t)=\xi_\alpha (t)-\xi_\alpha^{(N)}(t)$. Due to the ultrametric inequality, the jumps of $\eta_\alpha$ never exceed $p^N$ by absolute value, so that the process remains almost surely in $B_N$. It is proved in \cite{K2001} that the generator of the Markov process $\eta_\alpha$ on $B_N$ equals (on test function) $\DAN -\lambda I$ where
$$
\lambda =\frac{p-1}{p^{\alpha +1} -1}p^{\alpha (1-N)}.
$$
In \cite{K2001} (Theorem 4.9) the corresponding heat kernel is given explicitly.

In this paper we find an analytic interpretation of the latter operator using harmonic analysis on $B_N$ as an (additive) compact Abelian group (this group property, just as the above probabilistic construction, is of purely non-Archimedean nature and has no analogs in the classical theory of partial differential equations). We give an interpretation of $\DAN -\lambda I$ as a pseudo-differential operator on $B_N$, then consider it as an operator on $L^1(B_N)$ and study its Green function, the integral kernel of its resolvent. The choice of $L^1(B_N)$ as the basic space is motivated by applications to nonlinear equations.

The first model example of a nonlinear parabolic equation over $\Qp$ is the $p$-adic analog of the classical porous medium equation:
\begin{equation}
\frac{\partial u}{\partial t}+\DA (\Phi (u))=0,\quad u=u(t,x),\quad t>0,x\in \Qp,
\end{equation}
where $\Phi$ is a strictly monotone increasing continuous real function on $\mathbb R$. Its study was initiated in \cite{KK}. Here we consider this equation on $B_N$, taking the operator $\DAN$ instead of $\DA$:
\begin{equation}
\frac{\partial u}{\partial t}+\DAN (\Phi (u))=0.
\end{equation}

As in \cite{KK}, our study of Eq. (1.2) is based on general results by Crandall -- Pierre \cite{CP} and Br\'ezis -- Strauss \cite{BS} enabling us to consider this equation in the framework of nonlinear semigroups of operators. Following \cite{BV} we consider Eq. (1.2) also in $L^\gamma (B_N)$, $1<\gamma \le \infty$.

An important motivation of the present work is provided by the $p$-adic model of a poropus medium introduced in \cite{KOC1,KOC2}.

\medskip
\section{Preliminaries}

\medskip
{\it 2.1. p-Adic numbers} \cite{VVZ}.

Let $p$ be a prime
number. The field of $p$-adic numbers is the completion $\mathbb Q_p$ of the field $\mathbb Q$
of rational numbers, with respect to the absolute value $|x|_p$
defined by setting $|0|_p=0$,
$$
|x|_p=p^{-\nu }\ \mbox{if }x=p^\nu \frac{m}n,
$$
where $\nu ,m,n\in \mathbb Z$, and $m,n$ are prime to $p$. $\Qp$ is a locally compact topological field. By
Ostrowski's theorem there are no absolute values on $\mathbb Q$, which are not equivalent to the ``Euclidean'' one,
or one of $|\cdot |_p$.

The absolute value $|x|_p$, $x\in \mathbb Q_p$, has the following properties:
\begin{gather*}
|x|_p=0\ \mbox{if and only if }x=0;\\
|xy|_p=|x|_p\cdot |y|_p;\\
|x+y|_p\le \max (|x|_p,|y|_p).
\end{gather*}

The latter property called the ultra-metric inequality (or the non-Archi\-me\-dean property) implies the total disconnectedness of $\Qp$ in the topology
determined by the metric $|x-y|_p$, as well as many unusual geometric properties. Note also the following consequence of the ultra-metric inequality: $|x+y|_p=\max (|x|_p,|y|_p)$, if $|x|_p\ne |y|_p$.

The absolute value $|x|_p$ takes the discrete set of non-zero
values $p^N$, $N\in \mathbb Z$. If $|x|_p=p^N$, then $x$ admits a
(unique) canonical representation \index{Canonical
representation}
\begin{equation}
\label{2.1}
x=p^{-N}\left( x_0+x_1p+x_2p^2+\cdots \right) ,
\end{equation}
where $x_0,x_1,x_2,\ldots \in \{ 0,1,\ldots ,p-1\}$, $x_0\ne 0$.
The series converges in the topology of $\mathbb Q_p$. For
example,
$$
-1=(p-1)+(p-1)p+(p-1)p^2+\cdots ,\quad |-1|_p=1.
$$
We denote $\mathbb Z_p=\{ x\in \Qp:\ |x|_p\le 1\}$. $\mathbb Z_p$, as well as all balls in $\Qp$, is simultaneously open and closed.

Proceeding from the canonical representation (\ref{2.1}) of an element $x\in
\mathbb Q_p$, one can define the fractional part of $x$ as the rational number
$$
\{ x\}_p=\begin{cases}
0,& \text{if $N\le 0$ or $x=0$};\\
p^{-N}\left( x_0+x_1p+\cdots +x_{N-1}p^{N-1}\right) ,& \text{if
$N>0$}.\end{cases}
$$
The function $\chi (x)=\exp (2\pi i\{ x\}_p)$ is an additive
character of the field $\mathbb Q_p$, that is a character of its additive group. It is clear
that $\chi (x)=1$ if and only if $|x|_p\le 1$.

Denote by $dx$ the Haar measure on the
additive group of $\Qp $ normalized by the equality $\int_{\mathbb Z_p}dx=1$.

The additive group of $\Qp$ is self-dual, so that
the Fourier transform of a complex-valued function $f\in L^1(\Qp)$ is again a function on $\Qp$ defined as
$$
(\mathcal Ff)(\xi )=\int\limits_{\Qp}\chi (x\xi )f(x)\,dx.
$$
If $\mathcal Ff\in L^1(\Qp)$, then we have the inversion formula
$$
f(x)=\int\limits_{\Qp}\chi (-x\xi )\widetilde{f}(\xi )\,d\xi .
$$
It is possible to extend $\mathcal F$ from $L^1(\Qp)\cap L^2(\Qp)$ to a unitary operator on $L^2(\Qp)$, so that the Plancherel identity holds in this case.

In order to define distributions on $\Qp$, we have to specify a class of test functions. A function $f:\ \Qp\to \mathbb C$ is called locally constant if
there exists such an integer $l\ge 0$ that for any $x\in \Qp$
$$
f(x+x')=f(x)\quad \mbox{if }\|x'\|\le p^{-l}.
$$
The smallest number $l$ with this property is called the exponent of local constancy of the function $f$.

Typical examples of locally constant functions are additive characters, and also cutoff functions like
$$
\Omega (x)=\begin{cases}
1,& \text{if $\|x\|\le 1$};\\
0,& \text{if $\|x\|>1$}.\end{cases}
$$
In particular, $\Omega$ is continuous, which is an expression of the non-Archimedean properties of $\Qp$.

Denote by $\mathcal D(\Qp)$ the vector space of all locally constant functions with compact supports. Note that $\DD$ is dense in $L^q(\Qp)$ for each $q\in [1,\infty )$. In order to furnish $\DD$ with a topology, consider first the subspace $D_N^l\subset \DD$ consisting of functions with supports in a ball
$$
B_N=\{ x\in \Qp:\ |x|_p\le p^N\},\quad N\in \mathbb Z,
$$
and the exponents of local constancy $\le l$. This space is finite-dimensional and possesses a natural direct product topology. Then the topology in $\DD$ is defined as the double inductive limit topology, so that
$$
\DD=\varinjlim\limits_{N\to \infty}\varinjlim\limits_{l\to \infty}D_N^l.
$$

If $V\subset \Qp$ is an open set, the space $\mathcal D(V)$ of test functions on $V$ is defined as a subspace of $\DD$ consisting of functions with supports in $V$. For a ball $V=B_N$, we can identify $\mathcal D(B_N)$ with the set of all locally constant functions on $B_N$.

The space $\DDD$ of Bruhat-Schwartz distributions on $\Qp$ is defined as a strong conjugate space to $\DD$.

In contrast to the classical situation, the Fourier transform is a linear automorphism of the space $\DD$. By duality, $\FF$ is extended to a linear automorphism of $\DDD$. For a detailed theory of convolutions and direct products of distributions on $\Qp$ closely connected with the theory of their Fourier transforms see \cite{AKS,K2001,VVZ}.

\medskip
{\it 2.2. Vladimirov's operator} \cite{AKS,K2001,VVZ}.

The Vladimirov operator $\DA$, $\alpha >0$, of fractional differentiation, is defined first as a pseudo-differential operator with the symbol $|\xi|_p^\alpha$:
\begin{equation}
\label{2.2}
(\DA u)(x)=\FF^{-1}_{\xi \to x}\left[ |\xi |_p^{\alpha }\FF_{y\to \xi }u\right] ,\quad u\in \DD,
\end{equation}
where we show arguments of functions and their direct/inverse Fourier transforms. There is also a hypersingular integral representation giving the same result on $\DD$ but making sense on much wider classes of functions (for example, bounded locally constant functions):
\begin{equation}
\label{2.3}
\left( \DA u\right) (x)=\frac{1-p^\alpha }{1-p^{-\alpha -1}}\int\limits_{\Qp}|y|_p^{-\alpha -1}[u(x-y)-u(x)]\,dy.
\end{equation}

The Cauchy problem for the heat-like equation
$$
\frac{\partial u}{\partial t}+\DA u=0,\quad u(0,x)=\psi (x),\quad x\in\Qp,t>0,
$$
is a model example for the theory of $p$-adic parabolic equations. If $\psi$ is regular enough, for example, $\psi \in \DD$, then a classical solution is given by the formula
$$
u(t,x)=\int\limits_{\Qp}Z(t,x-\xi )\psi (\xi )\,d\xi
$$
where $Z$ is, for each $t$, a probability density and
$$
Z(t_1+t_2,x)=\int\limits_{\Qp}Z(t_1,x-y)Z(t_2,y)\,dy,\quad t_1,t_2>0,\ x\in \Qp.
$$

The "heat kernel" $Z$ can be written as the Fourier transform
\begin{equation}
Z(t,x)=\int\limits_{\Qp}\chi (\xi x)e^{-t|\xi |_p^\alpha}\,d\xi .
\end{equation}
See \cite{K2001} for various series representations and estimates of the kernel $Z$.

As it was mentioned in Introduction, the natural stochastic process in $B_N$ corresponds to the Cauchy problem
\begin{equation}
\frac{\partial u(t,x)}{\partial t}+\left( \DAN u\right) (t,x)-\lambda u(t,x)=0,\quad x\in B_N,t>0;
\end{equation}
\begin{equation}
u(0,x)=\psi (x),\quad x\in B_N,
\end{equation}
where the operator $\DAN$ is defined by restricting $\DA$ to functions $u_N$ supported in $B_N$ and considering the resulting function $\DA u_N$ only on $B_N$. Note that $\DAN$ defines a positive definite selfadjoint operator on $L^2(B_N)$, $\lambda$ is its smallest eigenvalue.

Under certain regularity assumptions, for example if $\psi \in \mathcal D(B_N)$, the problem (2.5)-(2.6) possesses a classical solution
$$
u(t,x)=\int\limits_{B_N}Z_N(t,x-y)\psi (y)\,dy,\quad t>0,x\in B_N,
$$
where
\begin{equation}
Z_N(t,x)=e^{\lambda t}Z(t,x)+c(t),
\end{equation}
$$
c(t)=p^{-N}-p^{-N}(1-p^{-1})e^{\lambda t}\sum \limits _{n=0}^\infty
\frac{(-1)^n}{n!}t^n\frac{p^{-N\alpha n}}{1-p^{-\alpha n-1}}.
$$
Another interpretation of the kernel $Z_N$ was given in \cite{CR}.

It was shown in \cite{KK} that the family of operators
$$
(T_N(t)u)(x)=\int\limits_{B_N}Z_N(t,x-y)\psi (y)\,dy
$$
is a strongly continuous contraction semigroup on $L^1(B_N)$. Its generator $A_N$ coincides with $\DAN -\lambda I$ at least on $\mathcal D(B_N)$. More generally, this is true in the distribution sense on restrictions to $B_N$ of functions from the domain of the generator of the semigroup on $L^1(\Qp )$ corresponding to $\DA$.

\medskip
\section{Harmonic analysis on the additive group of a $p$-adic ball}

Let us consider the $p$-adic ball $B_N$ as a compact subgroup of $\Qp$. As we know, any continuous additive character of $\Qp$ has the form $x\mapsto \chi (\xi x)$, $\xi \in \Qp$. The annihilator $\{ \xi \in \Qp :\ \chi (\xi x)=1$ for all $x\in B_N\}$ coincides with the ball $B_{-N}$. By the duality theorem (see, for example, \cite{Mo}, Theorem 27), the dual group $\widehat{B_N}$ to $B_N$ is isomorphic to the discrete group $\Qp /B_{-N}$ consisting of the cosets
\begin{equation}
p^m\left( r_0+r_1p+\cdots +r_{N-m-1}p^{N-m-1}\right) +B_{-N},\quad r_j\in \{ 0,1,\ldots ,p-1\},\quad m\in \mathbb Z,m<N.
\end{equation}
Analytically, this isomorphism means that any nontrivial continuous character of $B_N$ has the form $\chi (\xi x)$, $x\in B_N$, where $|\xi |_p>p^{-N}$ and $\xi \in \Qp$ is considered as a representative of the class $\xi +B_{-N}$. Note that $|\xi |_p$ does not depend on the choice of a representative of the class.

The normalized Haar measure on $B_N$ is $p^{-N}\,dx$. The normalization of the Haar measure on $\Qp /B_{-N}$ can be made in such a way (the normalized measure will be denoted $d\mu (x+B_{-N})$) that the equality
\begin{equation}
\int\limits_{\Qp}f(x)\,dx=\int\limits_{\Qp/B_{-N}}\left( p^N\int\limits_{B_{-N}}f(x+h)\,dh\right) d\mu (x+B_{-N})
\end{equation}
holds for any $f\in \mathcal D(\Qp )$; see \cite{Bo}, Chapter VII, Proposition 10; \cite{HR2}, (28.54). With this normalization, the Plancherel identity for the corresponding Fourier transform also holds; see \cite{HR2}, (31.46)(c).

On the other hand, the invariant measure on the discrete group $\Qp /B_{-N}$ equals the sum of $\delta$-measures concentrated on its elements multiplied by a coefficient $\beta$. In order to find $\beta$, it suffices to compute both sides of (3.2) for the case where $f$ is the indicator function of the set $\{ x\in \Qp :\ |x-p^{N-1}|_p\le p^{-N}\}$. Then the left-hand side equals $p^{-N}$ while the right-hand side equals $\beta$. Therefore $\beta =p^{-N}$.

The Fourier transform on $B_N$ is given by the formula
$$
(\FF_Nf)(\xi )=p^{-N}\int\limits_{B_N}\chi (x\xi )f(x)\,dx,\quad \xi \in (\Qp \setminus B_{-N})\cup \{0\},
$$
where the right-hand side, thus also $\FF_Nf$, can be understood as a function on $\Qp /B_{-N}$.

The fact that $\FF :\ \DD \to \DD$ implies that $\FF$ maps $\mathcal D(B_N)$ onto the set of functions on the discrete set $\widehat{B_N}$ having only a finite number of nonzero values. This set $\mathcal D(\widehat{B_N})$ with a natural locally convex topology can be seen as the set of test functions on $\widehat{B_N}=\Qp /B_{-N}$. The conjugate space $\mathcal D'(\widehat{B_N})$ consists of all functions on $\widehat{B_N}$ (see, for example, \cite{He}). Therefore the Fourier transform is extended, via duality, to the mapping from $\mathcal D'(B_N)$ to $\mathcal D'(\widehat{B_N})$. A theory of distributions on locally compact groups covering the case of $B_N$ was developed by Bruhat \cite{B}. To study deeper the operator $\DAN$, we need, within harmonic analysis on $B_N$, a construction similar to the well-known construction of a homogeneous distribution on $\Qp$ \cite{VVZ}.

Let us introduce the usual Riez kernel on $\Qp$,
$$
f_\alpha^{(N)}(x)=\frac{1-p^{-\alpha }}{1-p^{\alpha -1}}|x|_p^{\alpha -1},\quad \R \alpha >0,\quad \alpha \not \equiv 1 \pmod{\frac{2\pi i}{\log p}\mathbb Z}.
$$
Using the formula \cite{VVZ}
$$
\int\limits_{|x|_p\le p^N} |x|_p^{\alpha -1}dx=\frac{1-p^{-1}}{1-p^{-\alpha}}p^{\alpha N},
$$
we introduce a distribution from $\mathcal D'(B_N)$ setting
\begin{equation}
\left\langle
f_\alpha^{(N)},\varphi \right\rangle=\frac{1-p^{-1}}{1-p^{\alpha -1}}p^{\alpha N}\varphi (0)+\frac{1-p^{-\alpha}}{1-p^{\alpha -1}}\int\limits_{B_N}[\varphi (x)-\varphi (0)]|x|_p^{\alpha -1}dx,\quad \varphi \in \mathcal D(B_N).
\end{equation}

For $\R \alpha >0$, this gives
$$
\left\langle
f_\alpha^{(N)},\varphi \right\rangle= \frac{1-p^{-\alpha}}{1-p^{\alpha -1}}\int\limits_{B_N}|x|_p^{\alpha -1}\varphi (x)\,dx.
$$
On the other hand, the distribution (3.3) is holomorphic in $\alpha \not \equiv 1 \pmod{\frac{2\pi i}{\log p}\mathbb Z}$. Therefore  $f_{-\alpha}^{(N)}$ makes sense for any $\alpha >0$. Noticing that
$$
\frac{1-p^{-1}}{1-p^{-\alpha-1}}p^{-\alpha N}=\frac{p-1}{p^{\alpha +1}-1}p^{-\alpha N+\alpha}=\lambda
$$
(see Introduction), so that
\begin{equation}
\left\langle
f_{-\alpha}^{(N)},\varphi \right\rangle=\lambda \varphi (0)+\frac{1-p^{\alpha}}{1-p^{-\alpha -1}}\int\limits_{B_N}[\varphi (x)-\varphi (0)]|x|_p^{-\alpha -1}dx.
\end{equation}
The emergence of $\lambda$ in (3.4) ``explains'' its role in the probabilistic construction of a process on $B_N$ (\cite{K2001}, Theorem 4.9).

\medskip
\begin{teo}
The operator $\DAN$, $\alpha >0$, acts from $\mathcal D(B_N)$ to $\mathcal D(B_N)$ and admits, for each $\varphi \in \mathcal D(B_N)$, the representations:
\begin{itemize}
\item[\rm{(i)}]
$\DAN \varphi = f_{-\alpha}^{(N)}*\varphi$ where the convolution is understood in the sense of harmonic analysis on the additive group of $B_N$;
\item[\rm{(ii)}]
\begin{equation*}
\left( \DAN \varphi\right) (x)=\lambda \varphi (x)+\frac{1-p^{\alpha}}{1-p^{-\alpha -1}}\int\limits_{B_N}|y|_p^{-\alpha -1}[\varphi (x-y)-\varphi (x)]\,dy,\quad \alpha >0.
\end{equation*}
\item[\rm{(iii)}]
On $\mathcal D(B_N)$, $\DAN-\lambda I$ coincides with the pseudo-differential operator $\varphi \mapsto \FF_N^{-1}(P_{N,\alpha}\FF_N\varphi)$ where
\begin{equation}
P_{N,\alpha}(\xi )=\frac{1-p^{\alpha}}{1-p^{-\alpha -1}}\int\limits_{B_N}|y|_p^{-\alpha -1}[\chi (y\xi )-1]\,dy.
\end{equation}
This symbol is extended uniquely from $(\Qp \setminus B_{-N})\cup \{0\}$ onto $\Qp /B_{-N}$.
\end{itemize}
\end{teo}

\medskip
{\it Proof}. Denote, for brevity, $a_p=\dfrac{1-p^{\alpha}}{1-p^{-\alpha -1}}$. Let $x\in B_N$. Assuming that $\varphi$ is extended by zero onto $\Qp$, we find that
$$
\left( \DAN \varphi\right) (x)=a_p\int\limits_{\Qp}|y|_p^{-\alpha -1}[\varphi (x-y)-\varphi (x)]\,dy=I_1+I_2+I_3
$$
where
$$
I_1=a_p\int\limits_{B_N}|y|_p^{-\alpha -1}[\varphi (x-y)-\varphi (x)]\,dy,
$$
$$
I_2=a_p\int\limits_{|y|_p>p^N}|y|_p^{-\alpha -1}\varphi (x-y)\,dy,
$$
$$
I_3=-a_p\varphi (x)\int\limits_{|y|_p>p^N}|y|_p^{-\alpha -1}dy.
$$

We get using properties of $p$-adic integrals \cite{VVZ} that
$$
I_2=a_p\int\limits_{|x-z|_p>p^N}|x-z|_p^{-\alpha -1}\varphi (z)\,dz=a_p\int\limits_{|z|_p>p^N}|z|_p^{-\alpha -1}\varphi (z)\,dz=0;
$$
$$
I_3=-a_p\varphi (x)\sum\limits_{j=N+1}^\infty \int\limits_{|y|_p=p^j}|y|_p^{-\alpha -1}dy=-a_p\varphi (x)(1-\frac1p)\sum\limits_{j=N+1}^\infty p^{-\alpha j}=\lambda \varphi (x),
$$
which implies (ii). Comparing with (3.4) we prove (i).

In order to prove (3.5) we note that
\begin{multline*}
\FF_N\left( \DAN \varphi -\lambda \varphi \right)(\xi )=a_pp^{-N}\int\limits_{B_N}\chi (x\xi )\,dx\int\limits_{B_N}|y|_p^{-\alpha -1}[\varphi (x-y)-\varphi (x)]\,dy\\
=a_pp^{-N}\int\limits_{B_N}|y|_p^{-\alpha -1}dy\int\limits_{B_N}\chi (x\xi )[\varphi (x-y)-\varphi (x)]\,dx=P_{n,\alpha}(\xi )\left( \FF_N\varphi \right) (\xi ),
\end{multline*}
$\xi \in \Qp/B_{-N}. \qquad \blacksquare$

\bigskip
An important consequence of the representations given in Theorem 1 is the fact that, in contrast to operators on $\Qp$, $\DAN:\ \mathcal D(B_N)\to \mathcal D(B_N)$, so that we can define in a straightforward way, the action of this operator on distributions. In particular, the pseudo-differential representation remains valid on $\mathcal D'(B_N)$. Below (Theorem 3) this will be used to describe the domain of the operator $A_N$ on $L^1(B_N)$.

\medskip
\section{The Green function}

In Section 2 (just as in \cite{KK}) we defined the operator $A_N$ as the generator of the semigroup $T_N$ on $L^1(B_N)$. We can write its resolvent $(A_N+\mu I)^{-1}$, $\mu >0$, as
\begin{equation}
\left( (A_N+\mu I)^{-1}u\right) (x)=\int\limits_0^\infty e^{-\mu t}dt\int\limits_{B_N}Z_N(t,x-\xi )u(\xi )\,d\xi,\quad u\in L^1(B_N),
\end{equation}
where $Z_N$ is given in (2.7).

\medskip
\begin{teo}
The resolvent (4.1) admits the representation
\begin{equation}
\left( (A_N+\mu I)^{-1}u\right) (x)=\int\limits_{B_N}K_\mu (x-\xi )u(\xi )\,d\xi +\mu^{-1}p^{-N}\int\limits_{B_N}u(\xi )\,d\xi ,\quad u\in L^1(B_N),\mu >0,
\end{equation}
where for $0\ne x\in B_N$, $|x|_p=p^m$,
\begin{equation}
K_\mu (x)=\int\limits_{p^{-N+1}\le |\eta |_p\le p^{-m+1}}\frac{\chi (\eta x)}{|\eta |_p^\alpha -\lambda +\mu}\,d\eta .
\end{equation}
If $\alpha >1$, then for any $x\in B_N$,
\begin{equation}
K_\mu (x)=\int\limits_{|\eta |_p\ge p^{-N+1}}\frac{\chi (\eta x)}{|\eta |_p^\alpha -\lambda +\mu}\,d\eta .
\end{equation}
The kernel $K_\mu$ is continuous for $x\ne 0$ and belongs to $L^1(B_N)$.

If $\alpha >1$, then $K_\mu$ is continuous on $B_N$. If $\alpha =1$, then
\begin{equation}
|K_\mu (x)| \le C|\log |x|_p|,\quad x\in B_N.
\end{equation}
If $\alpha <1$, then
\begin{equation}
|K_\mu (x)| \le C|x|_p^{\alpha -1},\quad x\in B_N.
\end{equation}
\end{teo}

\medskip
{\it Proof}. Let us use the representation (2.7) substituting it into the equality
$$
\int\limits_{B_N}Z_N(t,x)\,dx=1
$$
(for the latter see Theorem 4.9 in \cite{K2001}). We find that
$$
c(t)=p^{-N}-e^{\lambda t}p^{-N}\int\limits_{B_N}Z(t,y)\,dy,
$$
so that
$$
Z_N(t,x)=e^{\lambda t}\left[ Z(t,x)-p^{-N}\int\limits_{B_N}Z(t,y)\,dy\right] +p^{-N},\quad x\in B_N.
$$

Let us consider the expression in brackets proceeding from the definition (2.4) of the kernel $Z$. Using the integration formula from Chapter 1, $\S$4 of \cite{VVZ} we obtain that
$$
Z(t,x)-p^{-N}\int\limits_{B_N}Z(t,y)\,dy=I_1(t,x)+I_2(t,x)
$$
where
$$
I_1(t,x)=\int\limits_{|\xi |_p\ge p^{-N+1}}\chi (\xi x)e^{-t|\xi |_p^\alpha}d\xi,
$$
$$
I_2(t,x)=\int\limits_{|\xi |_p\le p^{-N}}[\chi (\xi x)-1]e^{-t|\xi |_p^\alpha}d\xi,
$$
and $I_2(t,x)=0$ for $x\in B_N$.

Let $|x|_p=p^m$, $m\le N$. Then there exists such an element $\xi_0\in \Qp$, $|\xi_0|_p=p^{-m+1}$, that $\chi (\xi_0x)\ne 0$. Then making the change of variables $\xi =\eta +\xi_0$ we find using the ultra-metric property that
$$
\int\limits_{|\xi |_p\ge p^{-m+2}}\chi (x\xi )e^{-t|\xi |_p^\alpha}\,d\xi =\chi (x\xi_0 )\int\limits_{|\eta |_p\ge p^{-m+2}}\chi (x\eta )e^{-t|\eta |_p^\alpha}\,d\eta ,
$$
so that
$$
\int\limits_{|\xi |_p\ge p^{-m+2}}\chi (x\xi )e^{-t|\xi |_p^\alpha}\,d\xi =0.
$$
Therefore
$$
I_1(t,x)=\int\limits_{p^{-N+1}\le |\xi |_p\le p^{-m+1}}\chi (x\xi )e^{-t|\xi |_p^\alpha}\,d\xi ,
$$
thus
$$
Z_N(t,x)=e^{\lambda t}\int\limits_{p^{-N+1}\le |\xi |_p\le p^{-m+1}}\chi (x\xi )e^{-t|\xi |_p^\alpha}\,d\xi +p^{-N},\quad |x|_p=p^m.
$$

Substituting this in (4.1) and integrating in $t$ we come to (4.2) and (4.3). Note that $|\eta|_p^\alpha >\lambda$, as $|\eta|_p\ge p^{-N+1}$.

If $\alpha >1$, then the integral in (4.4) is convergent. For $|x|_p=p^m$ we prove repeating the above argument that
$$
\int\limits_{|\eta |_p\ge p^{-m+2}}\frac{\chi (\eta x)}{|\eta |_p^\alpha -\lambda +\mu}\,d\eta =0.
$$
Therefore in this case the representation (4.3) can be written in the form (4.4).

Obviously, $K_\mu (x)$ is continuous for $x\ne 0$. If $\alpha >1$, then there exists the limit
$$
\lim\limits_{x\to 0}K_\mu (x)= \int\limits_{|\eta |_p\ge p^{-N+1}}\frac{d\eta }{|\eta |_p^\alpha -\lambda +\mu}<\infty ,
$$
so that in this case $K_\mu$ is continuous on $B_N$.

Let $\alpha <1$. By (4.3) and an integration formula from \cite{VVZ}, Chapter 1, $\S$4,
\begin{multline*}
K_\mu (x)=\sum\limits_{l=-N+1}^{-m+1}\frac1{p^{\alpha l}-\lambda +\mu}\int\limits_{|\xi |_p=p^l}\chi (\xi x)\,d\xi \\
=(1-\dfrac1p)\sum\limits_{l=-N+1}^{-m}\frac{p^l}{p^{\alpha l}-\lambda +\mu}-
\frac{p^{-m}}{p^{\alpha (-m+1)}-\lambda +\mu},\quad |x|_p=p^m.
\end{multline*}

For some $\gamma >0$, $p^{\alpha l}-\lambda +\mu\ge \gamma p^{\alpha l}$. Computing the sum of a progression we obtain the estimate (4.6). Similarly, if $\alpha =1$, then $|K_\mu (x)|\le C(-m+N)$, which gives, as $m\to -\infty$, the inequality (4.5). $\qquad \blacksquare$

\bigskip
If $\alpha >1$, we can also give an interpretation of the resolvent $(A_N+\mu I)^{-1}$ in terms of the harmonic analysis on $B_N$. We have
\begin{equation}
(A_N+\mu I)^{-1}u=\left( K_\mu +\mu^{-1}\boldsymbol{1}\right) *u,\quad u\in L^1(B_N),
\end{equation}
where $\boldsymbol{1}(x)\equiv 1$, $K_\mu$ is given by (4.4), and the convolution is taken in the sense of the additive group of $B_N$.

Denote by $\Pi_N$ the set of all rational numbers of the form
$$
p^l\left( \nu_0+\nu_1p+\cdots +\nu_{-l+N-1}p^{-l+N-1}\right) ,\quad l<N,
$$
where $\nu_j\in \{ 0,1,\ldots ,p-1\}$, $\nu_0\ne 0$. As a set, the quotient group $\Qp/B_{-N}$ coincides with $\Pi_N\cup \{0\}$, and
$$
\{ \xi \in \Qp:\ |\xi |_p\ge p^{-N+1}\}=\bigcup\limits_{\eta \in \Pi_N}(\eta +B_{-N})
$$
where the sets $\eta +B_{-N}$ with different $\eta \in \Pi_N$ are disjoint.

Taking into account the fact that $\chi (\rho x)=1$ for $x\in B_N, \rho \in B_{-N}$, we find from (4.4) that
$$
K_\mu (x)=p^{-N}\sum\limits_{0\ne \eta \in \Qp/B_{-N}}\frac{\chi (\eta x)}{|\eta |_p^\alpha -\lambda +\mu}.
$$

Let us describe the domain $\Dom A_N$ of the generator of our semigroup $T_N(t)$ on $L^1(B_N)$ in terms of distributions on $B_N$.

\medskip
\begin{teo}
If $\alpha >1$, then the set $\Dom A_N$ consists of those and only those $u\in L^1(B_N)$, for which $f_{-\alpha}^{(N)}*u\in L^1(B_N)$ where the convolution is understood in the sense of the distribution space $\mathcal D'(B_N)$. If $u\in \Dom A_N$, then $A_Nu=f_{-\alpha}^{(N)}*u-\lambda u$ where the convolution is understood in the sense of the distributions from $\mathcal D'(B_N)$.
\end{teo}

\medskip
{\it Proof}. Let $u=(A_N+\mu I)^{-1}f$, $f\in L^1(B_N)$, $\mu >0$. Representing this resolvent as a pseudo-differential operator, we prove that $f_{-\alpha}^{(N)}*u-\lambda u+\mu u=f$ in the sense of $\mathcal D'(B_N)$.

Conversely, let $u\in L^1(B_N)$, $\DAN u=f_{-\alpha}^{(N)}*u\in L^1(B_N)$ where $\DAN$ is understood in the sense of $\mathcal D'(B_N)$. Set $f=(\DAN -\lambda I+\mu I)u$, $\mu >0$. Denote $u'=(A_N+\mu I)^{-1}f$. Then $u'\in \Dom A_N$, and the above argument shows that
$$
(\DAN -\lambda I+\mu I)(u-u')=0.
$$
Applying the pseudo-differential representation we see that
$$
\left[ P_{N,\alpha} (\xi )+\mu\right] \left[ (\FF_Nu)(\xi )-(\FF_Nu')(\xi )\right] =0,\quad \xi \in \Qp/B_{-N}.
$$
It is seen from (3.5) that the factor $P_{N,\alpha} (\xi )+\mu$ is real-valued, strictly positive and locally constant on $B_N$. Therefore the distribution $\FF_N u-\FF_N u'$ is zero. Since $\FF_N$ is an isomorphism (see \cite{B}), we find that $u=u'$, so that $u\in \Dom A_N$. $\qquad \blacksquare$

\medskip
\section{Nonlinear equations}

Let us consider the equation (1.2) where $\Phi$ is a strictly monotone increasing continuous real function, $\Phi (0)=0$, and the linear operator $\DAN$ is understood as the operator $A_N+\lambda I$ on $L^1(B_N)$. By the results from \cite{CP} and \cite{BS}, the nonlinear operator $\DAN \circ \Phi$ is $m$-accretive, which implies the unique mild solvability of the Cauchy problem for the equation (1.2) with the initial condition $u(0,x)=u_0(x)$, $u_0\in L^1(B_N)$; see e.g. \cite{Ba} for the definitions. As in the classical case \cite{BV}, this mild solution can be interpreted also as a weak solution.

Following \cite{BV}, we will show that the above construction of the $L^1$-mild solution gives also $L^\gamma$-solutions for $1<\gamma \le \infty$.

\medskip
\begin{teo}
Let $u(t,x)$, $t>0,x\in B_N$, be the above mild solution. If $0<u_0\in L^\gamma (B_N)$, $1\le \gamma \le \infty$, then $u(t,\cdot )\in L^\gamma (B_N)$ and
\begin{equation}
\|u(t,\cdot )\|_{L^\gamma (B_N)}\le \|u_0\|_{L^\gamma (B_N)}.
\end{equation}
\end{teo}

\medskip
{\it Proof}. The case $\gamma =1$ has been considered, while the case $\gamma =\infty$ will be implied by the inequality (5.1) for finite values of $\gamma$ (see Exercise 4.6 in \cite{Br}).

Thus, now we assume that $1<\gamma <\infty$. It is sufficient to prove (5.1) for $u_0\in \mathcal D(B_N)$. Indeed, if that is proved, we approximate in  $L^\gamma (B_N)$ an arbitrary function $u_0\in L^\gamma (B_N)$ by a sequence $u_{0,j}\in \mathcal D(B_N)$. For the corresponding solutions $u_j(t,x)$ we have
\begin{equation}
\|u_j(t,\cdot )\|_{L^\gamma (B_N)}\le \|u_{0,j}\|_{L^\gamma (B_N)}.
\end{equation}
Since our nonlinear semigroup consists of operators continuous on $L^1(B_N)$, we see that, for each $t\ge 0$, $u_j(t,\cdot )\to u(t,\cdot )$ in $L^1(B_N)$. By (5.2), the sequence $\{ u_j(t,\cdot )\}$ is bounded in $L^\gamma (B_N)$. These two properties imply the weak convergence $u_j(t,\cdot )\rightharpoonup u(t,\cdot )$ in $L^\gamma (B_N)$ (see Exercise 4.16 in \cite{Br}).

Next, we use the weak lower semicontinuity of the $L^\gamma$-norm (see Theorem 2.11 in \cite{LL}), that is the inequality
$$
\liminf\limits_j \|u_j(t,\cdot )\|_{L^\gamma (B_N)}\ge \|u(t,\cdot )\|_{L^\gamma (B_N)}.
$$
Passing to the lower limit in both sides of (5.2), we come to (5.1).

Let us prove (5.1) for $u_0\in \mathcal D(B_N)$, $1<\gamma <\infty$. By the Crandall-Liggett theorem (see \cite{Ba} or \cite{C}), $u(t,x)$ is obtained as a limit in $L^1(B_N)$,
$$
u(t,\cdot )=\lim\limits_{k\to \infty}\left( I+\frac{t}k\DAN \circ \Phi \right)^{-k}u_0,
$$
that is $u(t,\cdot )=\lim\limits_{k\to \infty}u_k$ where $u_k$ are found recursively from the relation
\begin{equation}
\frac{t}{k+1}\DAN \circ \Phi (u_{k+1})+u_{k+1}=u_k.
\end{equation}

Under our assumptions, $u(t,x)>0$ (this follows from Theorem 4 in \cite{CP}). The nonlinear operator $\left( I+\frac{t}k\DAN \circ \Phi \right)^{-1}$ is also positivity preserving (Proposition 1 in \cite{CP}), so that $u_k>0$ for all $k$.

Note that the operator $\DAN$ commutes with shifts while the equation (5.3) for $u_{k+1}$ has a unique solution in $L^1(B_N)$. As a result, if $u_0\in \mathcal D(B_N)$, then all the functions $u_k$ belong to $\mathcal D(B_N)$.

Rewriting (5.3) in the form
\begin{equation}
\left( \frac{t}{k+1}\right)^{-1}(u_{k+1}-u_k)=- \DAN \circ \Phi (u_{k+1}),
\end{equation}
multiplying both sides by $u_{k+1}^{\gamma -1}$ and integrating on $B_N$ we find that
\begin{equation}
\left( \frac{t}{k+1}\right)^{-1}\int\limits_{B_N}(u_{k+1}-u_k)u_{k+1}^{\gamma -1}dx=-\int\limits_{B_N} u_{k+1}^{\gamma -1}\DAN \circ \Phi (u_{k+1})\,dx.
\end{equation}

Let $w=u_{k+1}^{\gamma -1}$. Then $w\in \mathcal D(B_N)$. It follows from (5.4) that $\DAN \Phi (u_{k+1})\in \mathcal D(B_N)$. Also we have $\Phi (u_{k+1})\in \mathcal D(B_N)$, so that $\Phi (u_{k+1})$ belongs to the domain of a selfadjoint realization of the operator $\DAN$ in $L^2(B_N)$. Therefore we can transform the integral in the right-hand side of (5.5) as follows:
\begin{equation}
\int\limits_{B_N} u_{k+1}^{\gamma -1}\DAN \circ \Phi (u_{k+1})\,dx=\int\limits_{B_N}\Phi (w^{\frac1{\gamma -1}})\DAN (w)\,dx.
\end{equation}

The right-hand side of (5.6) is nonnegative by Lemma 2 of \cite{BS}. Now it follows from (5.5) that
$$
\int\limits_{B_N} u_{k+1}^\gamma dx\le \int\limits_{B_N} u_ku_{k+1}^{\gamma -1}dx.
$$
Applying the H\"older inequality we find that
$$
\int\limits_{B_N} u_{k+1}^\gamma dx\le \left( \int\limits_{B_N} u_k^\gamma dx\right)^{1/\gamma}\left( \int\limits_{B_N} u_{k+1}^\gamma dx\right)^{\frac{\gamma -1}\gamma},
$$
which implies the inequality
$$
\|u_{k+1}\|_{L^\gamma (B_N)}\le \|u_k\|_{L^\gamma (B_N)}
$$
and, by induction, the inequality
$$
\|u_{k+1}\|_{L^\gamma (B_N)}\le \|u_0\|_{L^\gamma (B_N)}.
$$
Passing to the limit, we prove (5.1). $\qquad \blacksquare$

\section*{Acknowledgment}
 This work was supported in part by Grant 23/16-18 ``Statistical dynamics, generalized Fokker-Planck equations, and their applications in the theory of complex systems'' of the Ministry of Education and Science of Ukraine.


\medskip
\begin{thebibliography}{999}
\bibitem{AKS}
S. Albeverio, A. Yu. Khrennikov and V. M. Shelkovich, {\it Theory of p-Adic Distributions. Linear and Nonlinear Models}. Cambridge University Press, 2010.
\bibitem{Ba}
V. Barbu, {\it Nonlinear Differential Equations of Monotone Types in Banach Spaces}, Springer, New York, 2010.
\bibitem{BV}
M. Bonforte and J. L. V\'azquez, Fractional nonlinear degenerate diffusion equations on bounded domains, {\it Nonlinear Anal.} {\bf 131} (2016), 363--398.
\bibitem{Bo}
N. Bourbaki, {\it Elements of Mathematics. Integration II}, Springer, Berlin, 2004.
\bibitem{Br}
H. Br\'ezis, {\it Functional Analysis, Sobolev Spaces and Partial Differential Equations}, Springer, New York, 2011.
\bibitem{BS}
H. Br\'ezis and W. Strauss, Semilinear elliptic equations in $L^1$, {\it J. Math. Soc. Japan} {\bf 25} (1973), 15--26.
\bibitem{B}
F. Bruhat, Distributions sur un groupe localement compact et applications \`a l'\'etude des repr\'esentations des groupes $p$-adiques, {\it Bull. Soc. Math. France} {\bf 89} (1961), 43--75.
\bibitem{CR}
O. F. Casas-S\'anchez and J. J. Rodr\'\i guez-Vega, Parabolic type equations on $p$-adic balls, {\it Boletin de Mat.} {\bf 22} (2015), 97--106.
\bibitem{C}
Ph. Cl\'ement et al., {\it One-Parameter Semigroups}, North-Holland, Amsterdam, 1987.
\bibitem{CP}
M. Crandall and M. Pierre, Regularizing effects for $u_t+A\psi (u)=0$ in $L^1$, {\it J. Funct. Anal.} {\bf 45} (1982), 194--212.
\bibitem{He}
A. Ya. Helemskii, {\it Lectures and Exercises on Functional Analysis}, AMS, Providence, 2006.
\bibitem{HR2}
E. Hewitt and K. A. Ross, {\it Abstract Harmonic Analysis, Vol. II}, Springer, Berlin, 1979.
\bibitem{KK}
A. Khrennikov and A. N. Kochubei, $p$-Adic analogue of the porous medium equation, {\it J. Fourier Anal. Appl.} (to appear), arXiv: 1611.08863.
\bibitem{KOC1}
A. Khrennikov, K. Oleschko and M. J. Correa Lopez, Application of p-adic wavelets to model reaction-diffusion dynamics in random porous media, {\it J. Fourier Anal. Appl.} {\bf 22} (2016), 809 -- 822.
\bibitem{KOC2}
A. Khrennikov, K. Oleschko and M. J. Correa Lopez, Modeling fluid's dynamics with master equations in ultrametric spaces representing the treelike
structure of capillary networks, {\it Entropy} {\bf 18} (2016), art. 249, 28 pp.
\bibitem{K2001}
A. N. Kochubei, {\it Pseudo-Differential Equations and Stochastics
over Non-Archimedean Fields}, Marcel Dekker, New York, 2001.
\bibitem{LL}
E. H. Lieb and M. Loss, {\it Analysis}, AMS, Providence, 2001.
\bibitem{Mo}
S. A. Morris, {\it Pontryagin Duality and the Structure of Locally Compact Abelian Groups}, Cambridge University Press, 1977.
\bibitem{VVZ}
V. S. Vladimirov, I. V. Volovich and E. I. Zelenov, {\it $p$-Adic Analysis and
Mathematical Physics}, World Scientific, Singapore, 1994.
\bibitem{V}
V. S. Vladimirov, {\it Tables of Integrals of Complex-Valued
Functions of $p$-Adic Arguments}, Steklov Mathematical Institute,
Moscow, 2003 (Russian). English version, ArXiv: math-ph/9911027.
\bibitem{Z2016}
W. A. Z\'u\~niga-Galindo, {\it Pseudodifferential Equations over Non-Archimedean Spaces}, Lect. Notes Math. Vol. 2174 (2016), XVI+175 p.

\end{thebibliography}
\end{document}